\documentclass[twoside,12pt,reqno]{amsart}
\usepackage{amsmath,amscd,amssymb,fullpage}
\usepackage{pstricks} 
\usepackage[all]{xy}
\newcommand{\wb}[1]{\overline{#1}}
\newcommand{\wt}{\widetilde}
\newcommand{\wh}{\widehat}
\newcommand{\Mod}{\operatorname{Mod}}
\newcommand{\Mor}{\operatorname{Mor}}
\newcommand{\End}{\operatorname{End}}
\newcommand{\Id}{\operatorname{Id}}
\newcommand{\cab}{\operatorname{sat}} 
\newcommand{\cut}{\operatorname{cut}}
\newcommand{\sch}{\operatorname{sch}}

\newcommand{\qd}{\operatorname{\mathsf{d}}}
\newcommand{\qdim}{\operatorname{qdim}}
\newcommand{\e}{\operatorname{e}}
\newcommand{\M}{{\mathcal{M}}}
\newcommand{\T}{{\mathcal{T}}}
\newcommand{\Tg}{{\mathcal{T}}_\Uh}
\newcommand{\ZZ}{\mathbb{Z}}  
\newcommand{\QQ}{\mathbb{Q}}
\newcommand{\RR}{\mathbb{R}}
\newcommand{\CC}{\mathbb{C}}
\newcommand{\NN}{\mathbb{N}}
\newcommand{\Sy}{\mathfrak{S}}
\newcommand{\go}{\longrightarrow}
\newcommand{\sll}{\mathfrak{sl}}
\newcommand{\slmn}{\sll(m|n)}
\newcommand{\Uh}{{{U_h}(\mathfrak{g})}}
\renewcommand{\H}{\mathcal H}
\newcommand{\qn}[1]{{\left\{#1\right\}}}
\newcommand{\qnc}[1]{{\left[#1\right]}}
\newcommand{\qum}[1]{\widetilde{#1}}
\newcommand{\h}{\ensuremath{\mathfrak{h}}}
\newcommand{\g}{\ensuremath{\mathfrak{g}}}
\newcommand{\D}{\mathbf D}

\newcommand{\wta}{\nu}
\newcommand{\Vschur}[1]{{{\mathsf S}^{#1}V}}
\newcommand{\mathsmall}[1]{\mbox{\small$#1$}} 

\newcommand{\Conway}{\nabla}
\newcommand{\fdeg}{\operatorname{f-deg}}

\newcommand{\cat}{{\mathcal{C}}}
\newcommand{\Fmn}{{{F}}_{m|n}}
\newcommand{\quo}[2]{\left<#1\,:\,#2\right>}

\newcounter{bibcount}

\newtheorem{prop}{\bf Proposition}[section]
\newtheorem{defi}[prop]{\bf Definition}
\newtheorem{lem}[prop]{\bf Lemma}
\newtheorem{theo}{\bf Theorem}
\newtheorem{coro}[prop]{\bf Corollary}
\newtheorem{conj}[prop]{\bf Conjecture}
\newtheorem{rk}[prop]{Remark}

\newcommand{\pspic}[2]{
  \psset{unit= #1 } 
  {\begin{array}{c} \hspace{-1.7ex}
      \raisebox{-2.25ex}{#2}
      \hspace{-0.85ex}\end{array}}}

\newcommand{\nbraid}[1]{
  \pspic{#1}{
    \begin{pspicture}(0,0)(1,1)
      \psline{<-}(.9,0)(.1,1)
      \psline{<-}(.1,0)(.4,.4)
      \psline{-}(.6,.6)(.9,1)
    \end{pspicture}}}
\newcommand{\pbraid}[1]{
  \pspic{#1}{
    \begin{pspicture}(0,0)(1,1)
      \psline{<-}(.1,0)(.9,1)
      \psline{<-}(.9,0)(.6,.4)
      \psline{-}(.4,.6)(.1,1)
    \end{pspicture}}}
\newcommand{\IIbraid}[1]{
  \pspic{#1}{
    \begin{pspicture}(0,0)(1,1)
      \pscurve{<-}(.1,0)(.4,.5)(.1,1)
      \pscurve{<-}(.9,0)(.6,.5)(.9,1)
    \end{pspicture}}}
\newcommand{\ptwist}[1]{
  \pspic{#1}{
    \begin{pspicture}(0,0)(1,1)
      \pscurve{<-}(.1,0)(.15,.45)(.45,.75)(.7,.5)(.4,.2)(.20,.38)
      \pscurve{-}(.14,.62)(.1,.9)(.1,1)
    \end{pspicture}}}
\newcommand{\Ibraid}[1]{
  \pspic{#1}{
    \begin{pspicture}(0,0)(1,1)
      \psline{<-}(.5,0)(.5,1)
    \end{pspicture}}}

\begin{document}
\title [Colored HOMFLY-PT and Multivariable Link Invariants] {On the
  Colored HOMFLY-PT, Multivariable and Kashaev Link Invariants} 
\author{Nathan Geer}
\address{School of Mathematics\\ 
Georgia Institute of Technology\\ 
Atlanta, GA 30332-0160}
\email{geer@math.gatech.edu}
\author{Bertrand Patureau-Mirand}
\address{LMAM, Universit\'e de Bretagne-Sud, BP 573\\
F-56017 Vannes, France }
\email{bertrand.patureau@univ-ubs.fr}
\date{\today}

\begin{abstract}
  We study various specializations of the colored HOMFLY-PT polynomial. These
  specializations are used to show that the multivariable link invariants
  arising from a complex family of $\slmn$ super-modules previously defined by
  the authors contains both the multivariable Alexander polynomial and
  Kashaev's invariants.  We conjecture these multivariable link invariants
  also specialize to the generalized multivariable Alexander invariants
  defined by Y.  Akutsu, T. Deguchi, and T. Ohtsuki.
\end{abstract}

\maketitle
\setcounter{tocdepth}{1}
\tableofcontents

\section*{Introduction}

Let $\sll(m|n)$ be the special linear Lie superalgebra.  In \cite{GP2} the
authors show that the Reshetikhin-Turaev quantum invariant of links arising
from the category of quantized $\sll(m,n)$ modules can be modified to produce
non-trivial multivariable link invariants for each $m,n\in\NN^*$ and $c\in
\NN^{m+n-2}$.  For $m\geq 2$, let $M_{\sll(m|1)}^0$ be the invariant
associated to $\sll(m|1)$ and ${0}\in\NN^{m-1}$.  If $L$ is a link with $k$
components and $k\geq 2$ then $M_{\sll(m|1)}^0(L)\in
\ZZ[q^{\pm1},q_1^{\pm1},\ldots,q_k^{\pm1}]$.  If $k=1$ then
$M_{\sll(m|1)}^0(L)\in g^{-1}\ZZ[q^{\pm1},q_1^{\pm1}]$ where $g$ is an element
of $\ZZ[q^{\pm1},q_1^{\pm1}]$ which does not depend on $L$.

The variable $q$ is the usual quantum parameter coming from the
Reshetikhin-Turaev quantum group construction.  The origin of the variable
$q_i$ is as follows.  The isomorphism classes of irreducible
finite-dimensional $\sll(m|1)$-module are parameterized by
$\NN^{m-1}\times\CC$.  The invariant $M_{\sll(m|1)}^0$ is constructed by
assigning the family of modules corresponding to $(0,...,0,\alpha_i)$, for
$\alpha_i \in \CC$, on the $i$th component of $L$.  The variable $q_i$ 
corresponds to $q^{\alpha_i}$ in this construction.

These invariants associate a variable to each component of the link.  There
are very few known invariants with such properties.  Such invariants including
the multivariable Alexander polynomial and the generalized multivariable
Alexander invariants $\{ADO_m\}_{m \geq 2}$ defined by Akutsu, Deguchi and
Ohtsuki \cite{ADO}.  In this paper we use various specializations of the
HOMFLY-PT polynomial to show that the invariants $\{M_{\sll(m|1)}^0\}_{m\geq
  2}$ contain the multivariable Alexander polynomial and Kashaev's invariants
\cite{Kv} (which are a specialization of the invariants $ADO_m$).

Using the quantum dilogarithm, Kashaev \cite{Kv} defines a family of complex
valued link invariants $K_{m}$ indexed by integers $m\geq 2$.  He then
observes that the hyperbolic volume of the complement of some simple
hyperbolic knots is determined by the asymptotic behavior of these link
invariants and conjectures that this is true for any hyperbolic knot.  In
\cite{MM}, H. Murakami and J. Murakami reformulate and strengthen Kashaev's
conjecture as follows: ``The colored Jones polynomials determine the
simplicial volume for any knot.''  This conjecture has become known as the
Volume Conjecture.

In making their reformulation H. Murakami and J. Murakami show that the set of
generalized multivariable Alexander invariants $\{ADO_{m}\}_{m\geq 2}$ and the
set of colored Jones polynomials have a non-trivial intersection.  Moreover,
they show that this intersection contains Kashaev's invariants.

In Section \ref{S:RLK} we will show that a similar result holds for the
invariants $M_{\sll(m|1)}^0$; namely, that the intersection of the set of
multivariable link invariants $\{M_{\sll(m|1)}^0\}_{m\geq 2}$ and the set of
colored HOMFLY-PT polynomials contains Kashaev's invariants.  In the first
half of this paper, we extend the fact that the two variable Links-Gould
invariants $\{LG^{m,1}\}_{m\geq 2}$ \cite{LG} specialize to the Alexander
polynomial \cite{WIL} by showing that the invariants $\{M_{\sll(m|1)}^0\}$
which can be seen as multivariable extensions of the Links-Gould invariants
specialize to the Conway function $\nabla$.

For a link $L$ with $k$ component and any integer $m\geq2$ the above discussion
can be summarized by the following diagram:
$$
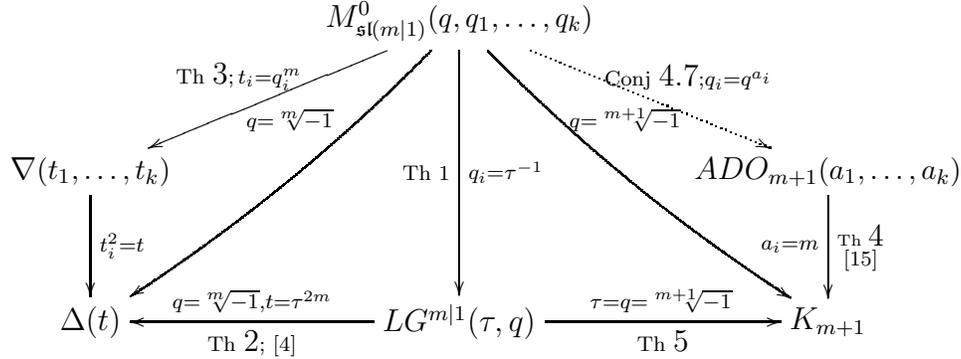
\xymatrix@!0 @R=11ex @C=27ex{ &
  M_{\sll(m|1)}^0\mathsmall{(q,q_1,\ldots,q_k)}
  \ar[ld]^-{q=\sqrt[m]{-1}}_-{\text{Th }\ref{th:M2Alex};\,t_i=q_i^m}
  \ar[dd]^-{q_i=\tau^{-1}}_-{\text{Th \ref{th:M2LG}}}
  \ar@{.>}[rd]_-{\,q=\sqrt[m+1]{-1}}
^-{\text{Conj }\ref{cj:M2ADO};
    q_i=q^{a_i}
}
  \ar@/^/[ddl]\ar@/_/[ddr]\\
\Conway\mathsmall{(t_1,\ldots,t_k)} \ar[d]^-{t_i^2=t}&&
ADO_{m+1}\mathsmall{(a_1,\ldots,a_k)} \ar[d]_-{a_i=m}^-{
\stackrel{\text{Th }\ref{th:ADO2K}}{\text{\cite{MM}}}}
\\
\Delta(t)& \ar[l]_-{q=\sqrt[m]{-1}, t=\tau^{2m}}^-{\text{Th
  }\ref{th:LG2Alex}\text{; \cite{WIL}}} LG^{m|1}(\tau,q)
\ar[r]^-{\tau=q=\sqrt[m+1]{-1}\,}_-{\text{Th }\ref{th:LG2K}}& K_{m+1} }$$
where a solid arrow represents an equality of the link invariants after a
specialization (and possibly a renormalization) and a dotted arrow represents
a similar equality which we conjecture to be true.  For example, the down
arrow in the middle of the diagram means
$$LG^{m,1}(\tau,q)=\left(\prod_{i=0}^{m-1}  \left(\tau^{-1} q^{i}-\tau
    q^{-i}\right) \right) M_m(q,\tau^{-1},\ldots, \tau^{-1})$$ and ``Th
\ref{th:M2LG}'' indicates that the proof of this equality is given in Theorem
\ref{th:M2LG}.

The construction of $\{M_{\sll(m|1)}^0\}$ differs from the construction of the
link invariants $\{ADO_m\}$.  The definition of $M_{\sll(m|1)}^0$ relies on
ribbon categories whereas the definition of $\{ADO_m\}$ use the Markov trace
for the colored braid group.  In both case, the standard method using ribbon
categories or the Markov trace is trivial.  To over come this difficulty, both
constructions rely on a regularization of the corresponding standard method.
The dotted arrow in the above diagram is correspond to Conjecture
\ref{cj:M2ADO}.

\section{The Multivariable Invariant associated with $\sll(m|1)$}\label{S:MI}
In this section we show that the invariant $ M_{\sll(m|1)}^0$ is related with
the Links-Gould invariant and the Conway function.

First, let us formulate a general notation which will be used throughout.  Let
$\cat$ be a category such that the set of endomorphisms of an object $V$ are a
module over some ring $K$.  Suppose that $f$ and $g$ are endomorphisms of $V$
such that $f=xg$ where $x\in K$ is a scalar.  Then we set $\quo f g=x$.

\subsection{$\Uh$-modules} 
Here we recall a particular category of modules over the Drinfeld-Jimbo
quantization associated to the special linear Lie superalgebra $\sll(m|n)$.
For more details please see \cite{G04A, GP2} and the references within.

Set $\g=\sll(m|n)$ ($m\neq n$) and let $\Uh$ be the DJ quantization associated
to $\g$ over $\CC[[h]]$.  In this paper, by a $\Uh$-module we mean a
topologically free $\Uh$-module of finite rank, i.e. a module over $\Uh$ which
is of the form $V[[h]]$ where $V$ is a finite dimensional $\g$-module.  In
\cite{G04A}, it is shown that every $\g$-module $V$ can be deformed to a
$\Uh$-module which we denote by $\qum{V}$.  Let $\M=\Mod_\Uh$ be the ribbon
category of topologically free $\Uh$-modules of finite rank and let $\Mod_\g$
be the category of finite dimensional $\g$-modules.
 
An object $\qum V$ of $\M$ is irreducible if $End_\Uh(\qum
V)=\CC[[h]]\Id_{\qum V}$.
 The classical limit of $\qum V$ is the $\g$-module $\qum
V/(h\qum V)$.  Denote by $C:\Mod_\Uh\go\Mod_\g$ the ``classical
limit'' functor.  Then $\qum V$ is a deformation of
$V=C(\qum V)$ (unique up to isomorphism).
\begin{lem}\label{lem:deform-irr}
  The deformation of an irreducible $\g$-module is an irreducible
  $\Uh$-module.
\end{lem}
\begin{proof}
  If $f\in\End_\Uh(\qum V)$, the weight decomposition of $V$ is preserved in
  $\qum V$ and $f$ respects it. So for a highest weight vector $v$ of $\qum
  V$, we have $f(v)=xv$ for some $x\in \CC[[h]]$.  Then $Ker(f-x\Id_{\qum V})$
  is a $\Uh$-module whose classical limit is $V$ implying $Ker(f-x\Id_{\qum
    V})=\qum V$ and so $f=x\Id_{\qum V}$.
\end{proof}
Lemma \ref{lem:deform-irr} implies that an irreducible $\Uh$-module is
determined by the underlying irreducible $\g$-module.  Every irreducible
finite-dimensional $\g$-module has a highest weight $\wta\in\h^*$ (where
$\h$ is the Cartan sub-superalgebra).  We denote such a module by $V_\wta$.
The set of isomorphism classes of irreducible finite-dimensional $\g$-modules
are in one to one correspondence with the set of dominant weights.  These
modules are parameterized by $\NN^{m+n-2}\times\CC$ and are divided into two
classes: typical and atypical.  We denote the weight corresponding to
$(c,\alpha)\in \NN^{m+n-2}\times\CC$ as $\wta_\alpha^c$ and say $\qum{V}$
is typical if $V$ is typical.
 
\subsection{The invariant $M^c_{\g}$} Let $\Tg$ be the category of
$\M$-colored ribbon graphs (that is the category of framed oriented
$\M$-colored graph with coupons labeled by even morphisms in $\M$).  Let
$\Fmn$ be the Reshetikhin-Turaev functor from $\Tg$ to $\M$ (see \cite{Tu94}).
In \cite{GP2}, the authors introduce a renormalization $\Fmn$ and show that
this renormalization leads to a multivariable invariant of ordered links
$M^c_{\g}$ for $c\in \NN^{m+n-2}$.  Let us now briefly recall this
construction.

Let $V$ be an irreducible module in $\M$.  Let $T_V\in\T_\M$ be a
$(1,1)$-ribbon graph whose open edge is oriented down and colored with $V$.
Then $\Fmn(T_V)\in \End_\M(V)=\CC[[h]]\Id_V$.  Using the notation at the
beginning of this section we have that $\Fmn(T_V)=\quo{\Fmn(T_V)}{\Id_V}\Id_V$.
  
For the braid closure $\hat T_V$ of $T_V$ (which is also its quantum trace in
$\T_\M$) one has $\Fmn(\hat T_V)=\qdim(V)\quo{\Fmn(T_V)}{\Id_V}$ which is zero
if $V$ is typical.  Hence $\Fmn$ is zero on any closed ribbon graph colored by
typical modules.
In \cite{GP2}, the authors define a modified quantum dimension which is a
function $\qd$ from the set of typical modules to $\CC[[h]][h^{-1}]$, which
has the property that the mapping $\hat T_V \mapsto
\qd(V)\quo{\Fmn(T_V)}{\Id_V}$ is a well defined invariant of $\M$-colored
closed ribbon graphs.  We denote this invariant by $F'$.

For a link $L=L_1\cup\cdots\cup L_k$ whose components are colored by the
deformations of typical modules $V_{\wta^c_{\alpha_i}}$ with the same
integral weight $c$ but with various complex parameters
$\alpha_1,\ldots\alpha_k$, $F'(L)$ depends continuously on these $k$
parameters.  In fact it defines an unique rational function
$M^c_{\g}(L)\in\QQ(q,q_1,\ldots,q_k)$ such that
\begin{equation}
  F'(L)=f(L).M^c_{\g}(L)(q,q^{\alpha_1},\ldots,q^{\alpha_k})
\end{equation}
where $f(L)$ is a function depending only on the linking matrix of $L$ (see
\cite{GP2}).  Here we use the notation
$$q=\e^{h/2}\qquad q^x=\e^{xh/2}\, .$$ 
When $L$ is not a knot (i.e. if $k\geq 2$),
$M^c_{\g}(L)$ happens to be a Laurent polynomial in these $k+1$
variables (for knots there is a denominator coming from $\qd$).
The function $M^c_{\g}$ is a multivariable invariant of (unframed)
ordered oriented links. 

For $\g=\sll(m|1)$, dominant weights are given by $\NN^{m-1}\times\CC$ which
can be written in the basis $(w_1,\ldots, w_m)$ of the fundamental weights.
As above irreducible finite-dimensional $\g$-modules are denoted
$V_{\wta^c_\alpha}$ ($c\in\NN^{m-1}$, $\alpha\in\CC$) where $
V_{\wta^0_\alpha}$ is typical iff $\alpha\notin\{0,-1,\ldots,1-m\}$.

Let $LG^{m|1}$ be the Links-Gould invariant (see \cite{WIL} and references
within).
\begin{theo}\label{th:M2LG}
  $$
  LG^{m|1}(q^{-\alpha},q)=\left(\prod_{i=0}^{m-1}
    \left(q^{\alpha+i}-q^{-(\alpha+i)}\right) \right)
  M^0_{\sll(m|1)}(q,q^\alpha,\ldots q^\alpha)
  $$
\end{theo}
\begin{proof}
  In \cite{WIL} the Links-Gould invariant is computed using an $R$-matrix of
  $U_q{\mathfrak{gl}}(m|1)\simeq U_q{\mathbb T}_1\otimes U_q\sll(m|1)$
  (isomorphism of Hopf algebras) where $U_q{\mathbb T}_1$ is the
  (co-)commutative Hopf algebra of polynomials in one primitive variable $t$.
  They consider the $2^{m}$-dimensional minimal typical
  $\sll(m|1)$-representation $\qum V_{\wta^0_\alpha}$ with highest weight
  $\alpha w_m=(0,\ldots,0,\alpha)$ (for a generic value of $\alpha$) on which
  $t$ acts by some scalar. Hence, they compute the Reshetikhin-Turaev
  invariant of a $(1,1)$-tangle $T_{\qum V_{\wta^0_\alpha}}$ where each
  component is colored by $\qum V_{\wta^0_\alpha}$ with a $R$-matrix that,
  differs from ours by a scalar and is rescaled so that the corresponding
  framed tangle invariant does not depend of its framing.  In other words, up
  to a correction for the framing the Links-Gould invariant of a link $L$ is
  given by $\quo{T_{\qum V_{\wta^0_\alpha}}}{\Id_{\qum
      V_{\wta^0_\alpha}}}$ where the closer of $T_{\qum
    V_{\wta^0_\alpha}}$ is the link $L$ whose components are all colored by
  $\qum V_{\wta^0_\alpha}$.
  
  Therefore, with the convention of \cite{WIL} where
  $\tau=q^{-\alpha}$, we get that for any link $L$,
  $$LG^{m|1}(L)(\tau,q)= 
  M^0_{\sll(m|1)}(L)(q,\tau^{-1},\ldots,\tau^{-1})/M^0_{\sll(m|1)}(unknot)
  (q,\tau^{-1}).$$ The theorem then follows from Lemma~\ref{L:d(V_a)1} and the
  fact that $M^0_{\sll(m|1)}(unknot) = \qd(\qum V_{\wta^0_\alpha})$.
\end{proof}

\begin{rk}
  Theorem \ref{th:M2LG} gives another proof that the Links-Gould invariant of
  a $(1,1)$-tangle $T$ depends only of the link closure $\wt T$ of $T$.  (This
  is not trivial for tangles with several components).
\end{rk}

In \cite{WIL} (Theorem 5) the following relation between $LG$ and the
Alexander-Conway polynomial is proved:
\begin{theo}\label{th:LG2Alex}(De Wit, Ashii, Links \cite{WIL}).  For all
  integer $m\geq2$, one has
  $$\Delta(L)(\tau^{2m})=LG^{m|1}(\tau,\e^{\sqrt{-1}\pi/m})$$
\end{theo}
We now prove a generalization of this theorem that shows that the invariants
$M^0_{\sll(m|1)}$ specialize to the multivariable Conway potential function.
For this, we use a modified version of Turaev's axioms for the Conway map
(\cite{Tu86} section 4) that the authors have proved and used in \cite{GP1}:
\begin{lem}\label{L:CAx}
  The Conway function is the map uniquely determined by
\begin{enumerate}
\item $\Conway$ assigns to each ordered oriented link $L$ in $S^3$ an element
  of the field $\QQ(t_1,\ldots,t_k)$ where $k$ is the number of components of
  $L$.
\item $\Conway(L)$ is unchanged under (ambient) isotopy of the link $L$.
\item $\Conway(\text{unknot})=(t_1-t_1^{-1})^{-1}$.
\item \label{Ax:int} If $k\geq2$ then
  $\Conway(L)\in\QQ[t_1^{\pm1},\ldots,t_k^{\pm1}]$.
\item The one variable function on links with several components
  $\wt\Conway(L)=\Conway(L)(t,t,\ldots,t)\in\QQ[t^{\pm1}]$ is unchanged by a
  reordering of the components of $L$.
\item (Conway Identity)
  $$
  \wt\Conway\left(\pbraid{4ex}\right) -\wt\Conway\left(\nbraid{4ex}\right)
  =(t-t^{-1})\wt\Conway\left(\IIbraid{4ex}\right)
  $$
\item \label{Ax:MDA}(Modified Doubling Axiom). If $L^+$ (resp. $L^-$) is
  obtained from the link $L=L_1\cup\cdots\cup L_k$ by replacing the
  $i$\textsuperscript{th} component $L_i$ by its $(2,1)$-cable (resp. by its
  $(2,-1)$-cable) then
  \begin{align*}
    t_k\Conway&(L^+)(t_1,\ldots,t_k)-t_i^{-1}\Conway(L^-)(t_1,\ldots,t_k)\,=\\
    &\left(\prod_{j\neq i} t_j^{lk_{ij}}\right)(t_i^2-t_i^{-2})
    \Conway(L)(t_1,\cdots,t_{i-1},t_i^2,t_{i+1},\cdots,t_k)
  \end{align*}
  where $(lk_{ij})_{i,j=1\cdots k}$ is the linking matrix of $L$.
\end{enumerate}
\end{lem}

\begin{theo}\label{th:M2Alex} For any oriented link $L$ with $k$ ordered
  component, one has
  $$\Conway(L)(q_1^m,\ldots,q_k^m)= \e^{\sqrt{-1}(m-1)\pi/2}
  M^0_{\sll(m|1)}(L)(\e^{\sqrt{-1}\pi/m},q_1,\ldots,q_k).$$
\end{theo}
\begin{proof}
  For $i\in\NN$, set $t_i=q_i^m$ and
  $M'(q_1,\ldots,q_k)=\e^{\sqrt{-1}(m-1)\pi/2}.
  M^0_{\sll(m|1)}(\e^{\sqrt{-1}\pi/m},q_1,\ldots,q_k)$ It is clear that $M'$
  satisfies the Axioms (1), (2), (4) and (5) of Lemma \ref{L:CAx} (we neglect
  the fact that $M'$ leave a priori in the extension $\QQ[q_i^{\pm1}]$ of
  $\QQ[t_i^{\pm1}]$).  Axiom (3) also follows easily from
  $$M'(q_1)(unknot)=
  \e^{\sqrt{-1}(m-1)\pi/2}. M^0_{\sll(m|1)}(\e^{\sqrt{-1}\pi/m},q_1)(unknot)=
  1/(q^m_1-q_1^{-m})$$ 
  where the second equality follows from Lemma \ref{L:d(V_a)}.

  Let us show that the Conway Identity and the Modified Doubling Axiom hold.
  To do this we need to recall the following facts.  For a link $L$ whose
  ordered components are colored by $\qum V_{\wta^0_{a_1}},\ldots,\qum
  V_{\wta^0_{a_k}}$, from \cite{GP2} we have that $F'$ and
  $M^0_{\sll(m|1)}$ are related by
  \begin{equation}
    \label{Eq:F'=M}
    \begin{array}{rcl}
      M^0_{\sll(m|1)}(L)(q,q^{a_1},\ldots,q^{a_k})&=&q^{-\sum lk_{ij}
        <a_iw_m,a_jw_m+2\rho>}F'(L)\\
      &=&q^{m\sum lk_{ij}a_i-r\sum lk_{ij}a_ia_j}F'(L)
    \end{array}
  \end{equation}
  where $w_m$ is the $m$\textsuperscript{th} fundamental weight, $<,>$ is the
  symmetric non-degenerate form on the Cartan sub-superalgebra defined in
  \cite{GP2}, and $r=<w_m,w_m>=\frac{m}{1-m}$.
  Let us temporarily extend the scalar to its quotient field. 
   
  Define $c_i\in \NN^{m-1}$ for $0\leq i \leq m$ as follows:
  $c_0=c_m=(0,...,0)$, $c_i=(0,...,0,1,0,...,0)$ where the $1$ is in the $i$th
  slot, for $1\leq i \leq m-1$.  One can use character formulas to show that
  (see \cite{GP2} (Lemma 2.8))~:
  \begin{equation}\label{Eq:split}
    \qum V_{\wta^0_\alpha}\otimes\qum V_{\wta^0_\alpha}\simeq
    \bigoplus_{i=0}^m\qum V_{i} 
  \end{equation}
  where $V_i$ is the irreducible module with highest weight $\wta_i$ which is
  equal to $\wta^{c_i}_{2\alpha+m-i-1}$ for $0\leq i\leq m-1$ and
  $\wta^{c_m}_{2\alpha}$ for $i=m$.  Thus, with the convention that $w_0=0$,
  we have $\wta_i=w_i+(2\alpha+m-i-1)w_m$.

  We consider a $2$-cable of the long Hopf link whose close component is
  colored by $\qum V_{\wta^0_\gamma}$ where $\gamma\in\CC$.  Its image under
  $F$ is an endomorphism $f$ of $\qum V_{\wta^0_\alpha}\otimes\qum
  V_{\wta^0_\alpha}$:
  $$
  f\ =\ F\left(\pspic{3ex}{
      \begin{pspicture}(0,0)(3.3,3)
        \psline{<-}(1,0)(1,.85) \psline{<-}(2,0)(2,.85)
        \pscurve{<-}(.85,2)(0,1.5)(1.5,1)(3,1.5)(2.15,2)
        \psline{-}(1.15,2)(1.85,2)
        \psline{-}(1,3)(1,1.15) \psline{-}(2,3)(2,1.15)
        \put(.4,.2){$\alpha$} \put(2.1,.2){$\alpha$} \put(3,1.6){$\gamma$}
      \end{pspicture}}\right).
  $$
  It acts on each simple summand $\qum V_i$ of $\qum
  V_{\wta^0_\alpha}\otimes\qum V_{\wta^0_\alpha}$ by a scalar $f_i$.  These
  scalars have been computed in \cite{GP2} (where they are denoted by $S'(\qum
  V_{\wta^0_\gamma},\qum V_i)$.  In particular,
  $f_i=\phi_{\wta_i+\rho}(\sch(V_{\wta^0_\gamma}))$ where $\sch$ is the super
  character of $V_{\wta^0_\gamma}$ and $\phi_{\wta_i+\rho}$ is a ring map.  In
  Proposition 1.1 of \cite{GP2} we show that
  $\sch(V_{\wta^0_\gamma})=u^{\gamma w_m}\chi'_1$ where the coefficient of
  $u^a$ is the dimension of the $a$-weight space and $\chi'_1$ is a linear
  combination of elements $u^b$ which are indexed by the set of positive odd
  roots $\{b\}$.  It follows that $f_i= q^{2\gamma(i-m)}q^{r(2\gamma+m-1)}
  \phi_{\wta_i+\rho}(\chi'_1)$ and $\phi_{\wta_i+\rho}(\chi'_1)$ is a non-zero
  Laurent polynomial in two variables evaluated at $q$ and $q^{2\alpha}$.
  
  Hence in a ring where the $f_i-f_j$ are invertible, (for example, for
  irrational values of $\alpha,\gamma\in\CC$ when $q=\e^{\sqrt{-1}\pi/m}$),
  the projector $P_i$ on $\qum V_i\simeq Im(P_i)$ is realized in $\End(\qum
  V_{\wta^0_\alpha}\otimes\qum V_{\wta^0_\alpha})$ as linear combination of
  power of $f$.  Hence the decomposition of Equation \eqref{Eq:split} is still
  valid for $q=\e^{\sqrt{-1}\pi/m}$.
  
  Now any endomorphism $g$ of $\qum V_{\wta^0_\alpha}\otimes\qum
  V_{\wta^0_\alpha}$ acts on $\qum V_i$ by some scalar $g_i$.  Hence $g=\sum_i
  g_iP_i$.  For the braid closure $\wh{P_i}$ of $P_i$, one has
  $F'(\wh{P_i})=\qd(\qum V_i)$ and so $F'(\wh g)=\sum_ig_i\qd(\qum V_i)$.  For
  $q=\e^{\sqrt{-1}\pi/m}$ and $0<j<m$, one has $\qd(\qum V_j)=0$ (see Lemma
  \ref{L:d(V_i)}).  So modulo the kernel of $F'$, one has, $\dim_\CC(\End(\qum
  V_{\wta^0_\alpha}\otimes\qum V_{\wta^0_\alpha}))=2$.  Using this, we are
  going to show that for $q=\e^{\sqrt{-1}\pi/m}$, $F'$ satisfies the following
  two skein relations:
 \begin{equation}
    \label{Eq:conway}
    q^{m\alpha-r\alpha^2}F'\left(\pbraid{4ex}\right)
    -q^{-m\alpha+r\alpha^2}F'\left(\nbraid{4ex}\right)= 
    (q^{m\alpha}-q^{-m\alpha})F'\left(\IIbraid{4ex}\right)
  \end{equation}
  \begin{equation}
    \label{Eq:cable}
    q^{2m\alpha-r\alpha^2}F'\left(\pbraid{4ex}\right)
    -q^{-2m\alpha+r\alpha^2}F'\left(\nbraid{4ex}\right)=
    (q^{2m\alpha}-q^{-2m\alpha})P_m
  \end{equation}
  where $r=<w_m,w_m>=\frac{m}{1-m}$ as above and all the components of the
  tangles are colored by $\qum V_{\wta^0_\alpha}$.
  
  To compute the coefficients in these two relations, one can consider the
  highest weight vector $v_+$ and the lowest weight vector $v_-$ of $\qum
  V_{\wta^0_\alpha}$.  Their weight are respectively $\alpha w_m$ and
  $(\alpha+m-1)w_m$.  The vector $v_+$ is even but the parity of $v_-$ is the
  parity of $m$. From the character formula one can also see that $v_+\otimes
  v_+$ is a highest weight vector of $\qum V_m$ and $v_-\otimes v_-$ is a
  lowest weight vector of $\qum V_0$.  So the two skein relations can be
  checked by evaluating the corresponding morphisms of $\qum
  V_{\wta^0_\alpha}\otimes\qum V_{\wta^0_\alpha}$ on $v_+\otimes v_+$ and
  $v_-\otimes v_-$.  Now the quasi-R-matrix acts by $1$ on these two vectors,
  so the positive braiding is simply given by
  $$
  F\left(\pbraid{4ex}\right)(v_+\otimes v_+)=q^{<\alpha w_m,\alpha
    w_m>}v_+\otimes v_+
  $$
  $$
  F\left(\pbraid{4ex}\right)(v_-\otimes v_-)=
  (-1)^mq^{<(\alpha+m-1)w_m,(\alpha+m-1)w_m>}v_-\otimes v_-
  $$
  where the $(-1)^m$ sign comes from the parity of $v_-$.  The negative
  braiding acts by the inverse values and the skein relations
  \eqref{Eq:conway} and \eqref{Eq:cable} are consequences of this.
  \\

  The skein relation \eqref{Eq:conway} implies that $M'$ satisfies the axiom
  of the Conway identity.
  \\

  We now use \eqref{Eq:cable} to check the modified doubled axiom for the
  oriented link $L$ with its $k$ ordered components colored by $\qum
  V_{\wta^0_{\alpha_i}}$, $i=1\cdots k$.  We choose a framing on $L$ so that
  $lk_{ii}=0$.  Hence $L^{\pm}$ can be obtain from $L$ by replacing $L_i$ by
  two parallel copies modified in a small ball as in the left hand side of
  \eqref{Eq:cable}.  This gives
  $$
  q^{2m\alpha_i-r\alpha_i^2}F'(L^+)-q^{-2m\alpha_i+r\alpha_i^2}F'(L^-)=
  (q^{2m\alpha_i}-q^{-2m\alpha_i})F'(\wb L)
  $$
  where the $j$-th component of $L^+$, $L^-$ and $\wb L$ is colored by
  $\alpha_j$ except the $i$-th component of $\wb L$ which is colored by
  $2\alpha_i$. Now as $lk_{ij}(L^\pm)=lk_{ji}(L^\pm)=2lk_{ij}(L)$ for $j\neq
  i$ and $lk_{ii}(L^\pm)=\pm1$, the framing correction of \eqref{Eq:F'=M}
  gives
  \begin{align*}
    t^{\alpha_i}M'(L^+)&(t^{\alpha_1},\ldots,t^{\alpha_k})
    -t^{-\alpha_i}M'(L^-)(t^{\alpha_1},\ldots,t^{\alpha_k})=\\
    &\left(\prod_{j\neq i}
      t^{lk_{ij}{\alpha_j}}\right)(t^{2\alpha_i}-t^{-2\alpha_i})
    M'(L)(t^{\alpha_1},\cdots,t^{\alpha_{i-1}},t^{2\alpha_i},t^{\alpha_{i+1}},
    \cdots,t^{\alpha_k})
  \end{align*}
  where $t=q^m$.  Hence $M'$ is the Conway function.
\end{proof}
\section{The Hecke Category $\H$}
In this section we introduce the Hecke category and consider its relations with
the ribbon category of $U_h\sll(m|n)$-module. 
We fix an integer $N\geq3$.  
\subsection{Definition of the Hecke category}
Consider the quotient field $K$ of the ring of polynomials in three variable
$a$, $s$ and $v$. Let $\qnc{k}_s=\dfrac{s^k-s^{-k}}{s-s^{-1}}$ for $k\in\ZZ$
and $\qnc{k}_s!=\qnc{k}_s\qnc{k-1}_s\cdots\qnc{1}_s$ for $k\in\NN$.
Let $R$ be the sub-ring of $K$ defined by
$$R=\QQ\left[a^{\pm1},s^{\pm1},v^{\pm1},
\dfrac{v-v^{-1}}{s-s^{-1}},(\qnc{N-1}_s!)^{-1}\right]$$ 

We first briefly recall the definition of the $R$-linear category of framed
oriented tangles $\T$.  Let $\D^2=\{z\in\CC:\,|z|\leq1\}$ be the standard disc
in $\CC\simeq\RR^2$ and let $M$ be the free associative monoid in two symbol
$\{+;-\}$.  To $(\alpha,\beta)\in M^2$ we associate the oriented
$0$-dimensional submanifold $\{(\frac0p,0),\ldots,(\frac{p-1}p,0)\}\cup
\{(\frac0q,1),\ldots,(\frac{q-1}q,1)\}\subset\D^2\times [0,1]$ where $p$ and
$q$ are the respective length of $\alpha$ and $\beta$. This set of points is
oriented according to $\alpha$ and the opposite of $\beta$.

The set of object of $\T$ is $M$.  The set of morphisms from $\alpha$ to
$\beta$ is the $R$-linear space span by isotopy class of framed oriented
tangles $T$ in $\D^2\times [0,1]$ such that the orientation of the boundary
$\partial T=T\cap\D^2\times \{0;1\}$ is given by $(\alpha,\beta)$ and the
framing on a point of $\partial T$ is given by the vector $(i,0)$.  We denote
this set of morphism by $\T(\alpha,\beta)$.  The composition of
$T_1\in\T(\alpha,\beta)$ with $T_2\in\T(\beta,\gamma)$ is obtained by gluing
the pairs $(\D^2\times [0,1],T_1)$ and $(\D^2\times [0,1],T_2)$ where
$(\D^2\times \{0\},\beta)$ is identified with $-(\D^2\times \{1\},-\beta)$.
The tensor product for objects comes from the monoidal structure of $M$ and
the tensor product of morphisms is induced by trivial embeddings $\D^2\times
[0,1]\amalg\D^2\times [0,1]\hookrightarrow\D^2\times [0,1]$.

We follow \cite{CB} to define a framed version of the Hecke category which we
denote by $\H$~: $\H$ is the quotient of $\T$ by the HOMFLY-PT skeins
relations \psset{unit= 6ex }
\begin{center}
  \vspace{1ex}
  $a^{-1}\pbraid{6ex}\ -\ a\nbraid{6ex}\ =\ (s-s^{-1})\IIbraid{6ex}$
  \\[2ex]
  $\ptwist{6ex}\ =\ av^{-1}\Ibraid{6ex}\hspace{12ex} L\ \cup\
  \pspic{6ex}{\begin{pspicture}(0,0)(1,1) \pscircle(.5,.5){.4}
      \psline{->}(.12,.5)(.12,.47)
  \end{pspicture}}\ =\ {\frac{v^{-1}-v}{s-s^{-1}}}\ L$.
\vspace{1ex}
\end{center}
We denote $F_\H$ the quotient functor from $\T$ to $\H$.

\subsection{Relation with $U_q\sll(m|n)$}
Let $V$ be the standard representation of $\g=\sll(m|n)$ of dimension
$m+n$ (which is an atypical irreducible module).
It is a well known fact that $\Fmn$ restricted to tangles whose components are
all colored by $\qum V$ satisfies a HOMFLY-PT-type skein relation.  More
precisely, let $F_V:\T\go\Tg$ be the functor that colors each component of a
tangle with $\qum V$.  For $\delta\in\ZZ^*$, let
$\psi_\delta:R\go\QQ[q^{\pm1/\delta},\qnc{N-1}_q!^{-1}]\subset\QQ[[h]]$ be the
ring morphism defined by
$$\psi_\delta(s)=q,\quad \psi_\delta(v)=q^{-\delta}\quad
\text{and}\quad\psi_\delta(a)=q^{-1/\delta}.$$ 
\begin{prop}\label{P:H2sl}
  There exists an unique monoidal functor $\wb \Fmn:\H\go\Mod_\Uh$ such that
  the following diagram is commutative
  $$\xymatrix @R=8ex @C=12ex{
    \ar[d]^{F_\H}\T\ar[r]^{F_V}&\ar[d]^{\Fmn}\Tg\\
    \H\ar[r]^{\wb \Fmn}&\Mod_\Uh}.$$ 
  Furthermore, $\wb \Fmn$ is $R$-linear in the following way:
  $$\forall x\in R,\,\forall T\in\Mor(\H),\, \wb
  \Fmn(xT)=\psi_{m-n}(x) \wb \Fmn(T).$$
\end{prop}

\subsection{Framing and grading}
There is a little difference between the framed and unframed version of
HOMFLY-PT.  One recovers the unframed skein relations by taking $a=v$.  But
this quotient makes the statement of Proposition \ref{P:H2sl} and the cabling
process less natural.  So we prefer to use a framed version and introduce a
``framing-degree'' ($\fdeg$) to make the correspondence more clear.

The category $\T$ possess a natural $\ZZ$-grading: if $T$ is a tangle define
$\fdeg(T)\in \ZZ$ as the total algebraic number of crossing of $T$.  For the
ring elements, we consider the $\ZZ$-grading given by $\fdeg(a^{\pm1})=\pm1$,
$\fdeg(v^{\pm1})=\fdeg(s^{\pm1})=0$.  Thus the spaces of morphism of $\T$ are
$\ZZ$-graded vector spaces.  One can easily check that the HOMFLY-PT relations
are homogeneous and thus $\H$ inherits this $\ZZ$-grading.  If $L$ is a framed
link with $f=\fdeg(L)$ then $F_\H(L)$ is an element of
$a^f\ZZ\left[s^{\pm1},v^{\pm1},\frac{v-v^{-1}}{s-s^{-1}}\right]$.

\section{The colored HOMFLY-PT polynomials dominate $M^c_{\sll(m|n)}$}
\subsection{Idempotents of the Hecke algebra}\label{S:Hecke-idempot}
Let $r\in\NN^*$.  The algebra $\H(+^r,+^r)$ is isomorphic to the Hecke algebra
$H_r$ of type $A$.  The simple subfactors of $H_r\otimes K$ are indexed by
Young diagrams $\lambda$ of size $|\lambda|=r$.  In other words, one has
$H_r\subset H_r\otimes K\simeq\bigoplus_{|\lambda|=r} H_{\lambda}$ where
$H_{\lambda}$ is isomorphic to the ring of $d(\lambda)$ by $d(\lambda)$ matrices
over $K$ and $d(\lambda)$ is the integer $r!$ divided by the product of the
hook lengths of cells of $\lambda$
(we say that a cell $c=(i,j)\in\NN^*\times\NN^*$ is in $\lambda$ if
$j\leq\lambda_i$. Then the hook length of $c$ is
$hl(c)=\lambda_i+\lambda'_j-i-j+1$ with $\lambda'$ the conjugate partition of
$\lambda$).

Following \cite{CB}, we use the notation $\qnc{hl(\lambda)}_s$ for the product
over all cells $c$ of $\lambda$ of their quantum hook-lengths:
$\qnc{hl(\lambda)}_s=\prod_ {c\in\lambda}\qnc{hl(c)}_s$.  We say that
$\lambda$ is $R$-admissible if $\qnc{hl(\lambda)}_s$ is invertible in $R$.
Then the projector on the simple subfactor of type $\lambda$ is realized by
multiplication with a central idempotent of $c_{\lambda}\in H_r$
($c_{\lambda}$ is the sum of the $d(\lambda)$ minimal orthogonal idempotents
of type $\lambda$ constructed in \cite{AM,CB}).  Moreover, $c_{\lambda} H_r$
is isomorphic to the ring of $d(\lambda)$ by $d(\lambda)$ matrices over $R$.
For any $R$-admissible diagram $\lambda$, we choose a minimal idempotent
$y_{\lambda}$ of type $\lambda$ (i.e. $y_{\lambda}\in c_{\lambda} H_r$).  This
choice is unique up to conjugation and $\fdeg(y_{\lambda})=0$ (see
\cite{AM,CB}).

\subsection{Classical limit}
Fix $\delta\in\ZZ^*$. We consider the surjective morphism of
$\QQ$-algebras $\pi$ from $H_r$ to the group ring of the permutation group
$\QQ[\Sy(r)]$ given by:
$$\pi(a)=1,\quad \pi(v)=1,\quad \pi(s)=1,\quad
\pi\left(\dfrac{v-v^{-1}}{s-s^{-1}}\right)=-\delta.$$ 
Then $\pi(y_{\lambda})$ is a minimal idempotent of the simple factor of
$\QQ[\Sy(r)]$ of type $\lambda$.
We now use the super version of the Schur functors ``$\Vschur \lambda$'' in
the following proposition which is the quantum analog of a theorem of Sergeev
(\cite{Se} Theorem 2):
\begin{prop}\label{P:part2weight}
  $\wb \Fmn(y_{\lambda})$ is a minimal projection on an irreducible
  $\Uh$-module $\qum {\Vschur{\lambda}}\subset\qum V^{\otimes r}$.  Write
  $\lambda=(\lambda_1,\lambda_2,\lambda_3,\ldots)$ (with
  $\lambda_i\leq\lambda_{i+1}$ and $\lambda_i=0$ for $i>>0$) and denote by
  $\lambda'=(\lambda'_1,\lambda'_2,\lambda'_3,\ldots)$ the conjugate partition
  of $\lambda$ and $\lambda''_i=\max(\lambda'_i-m,0)$.  Then,
\begin{enumerate}
\item if $\lambda_{m+1}>n$ then $\qum{\Vschur{\lambda}}={0}$.
\item If $\lambda_{m+1}\leq n$ then $\qum{\Vschur{\lambda}}$ is a deformation
  of an irreducible $\g$-module with highest weight $(\lambda_1-\lambda_2,
  \lambda_2-\lambda_3, \ldots, \lambda_{m-1}-\lambda_m, \lambda_m-\lambda''_1,
  \lambda''_1-\lambda''_2, \lambda''_2-\lambda''_3, \ldots,
  \lambda''_{n-1}-\lambda''_n)$.
\item In particular if $\lambda_{m+1}\leq n$ and $\lambda_{m}\geq n$
  then $\qum V^\lambda$ has highest weight $(\lambda_1-\lambda_2,
  \lambda_2-\lambda_3, \ldots, \lambda_{m-1}-\lambda_m,
  \lambda_m-\lambda'_1+m, \lambda'_1-\lambda'_2,
  \lambda'_2-\lambda'_3, \ldots, \lambda'_{n-1}-\lambda'_n)$
\end{enumerate}
\end{prop}
\begin{proof}
  Recall that $C$ the functor ``classical limit'' is obtained by sending $h$
  to zero:
  $$C: \Mod_\Uh\go\Mod_\g \;\text{ given by } \; \qum V\mapsto\qum
    V/(h\qum V).$$ 
  In the following, we will identify the $\g$-module $V$ with $C(\qum V)$.
  The representation of $H_r\simeq \H(+^r,+^r)$ given by $C\circ \wb \Fmn$
  factor as $\rho\circ \pi$ where $\pi:H_r\go\QQ[\Sy(r)]$ is the projection of
  $H_r$ onto the group ring of the permutation group $\Sy(r)$ and $\rho$ is
  the super representation of $\Sy(r)$ in $V^{\otimes r}$ (this is true
  because the universal $R$-matrix of $\Uh$ is $1$ modulo $h$).

  Since $\pi(y_{\lambda})$ is a minimal idempotent of the simple factor of
  $\QQ[\Sy(r)]$ of type $\lambda$, $\rho\circ \pi(y_{\lambda})$ is a
  projection on an irreducible $\g$-submodule $\Vschur{\lambda}$ of
  $V^{\otimes r}$ (see \cite{Se}). Hence $\wb \Fmn(y_{\lambda})$ is a
  projection on a $\Uh$-module $\qum{\Vschur{\lambda}}\subset \qum{V}^{\otimes
    r}$ whose classical limit is $\Vschur{\lambda}$. Since $\Vschur{\lambda}$
  is irreducible it follows that $\qum{\Vschur{\lambda}}$ is irreducible.  The
  theorem then follows from the description of the $\mathfrak{gl}(m|n)$-module
  $\Vschur{\lambda}$ made by Sergeev (\cite{Se} Theorem 2).
\end{proof}
If $L=(L_1\cup\cdots\cup L_k)$ is a framed oriented link with $k$ components,
and $\lambda^*=(\lambda^1,\ldots,\lambda^k)$ is a $k$-tuple of $R$-admissible
Young diagrams, one can define the $\lambda^*$-satellite of $L$, denoted
$\cab_{\lambda^*}(L)$, as the linear combination of links obtained by
replacing a regular neighborhood of $L_i$ by a torus that contains the braid
closure of $y_{\lambda^i}$, for each $1\leq i \leq k$.
Suppose now that $L$ is the closure of the (unique
up to isotopy) $(1,1)$-tangle $\cut(L,L_i)\in\T(+,+)$ obtained by opening the
component $L_i$.  Then $\cut(L,L_i)$ inherits the coloration $\lambda^*$ and
one can define the $\lambda^*$-cable of $\cut(L,L_i)$ as
$\cab_{\lambda^*}(\cut(L,L_i))\in\T(+^{|\lambda^i|},+^{|\lambda^i|})$.

A simple count shows that $\fdeg(\cut(L,L_i))=\fdeg(L)$ and
$$\fdeg(\cab_{\lambda^*}(\cut(L,L_i)))=\fdeg(\cab_{\lambda^*}(L))=
{^t|\lambda^*|}.lk(L).|\lambda^*|$$
where $lk(L)$ is the linking matrix of $L$ and $|\lambda^*|$ is the vector
column $(|\lambda^1|,|\lambda^2|,\ldots,|\lambda^k|)$.
As $y_{\lambda^i}$ is an idempotent in $\H$,
$$F_\H(\cab_{\lambda^*}(\cut(L,L_i)))\circ
y_{\lambda^i}=y_{\lambda^i}\circ F_\H(\cab_{\lambda^*}(\cut(L,L_i)))
=F_\H(\cab_{\lambda^*}(\cut(L,L_i))).$$ 
Now since $y_{\lambda^i}$ is minimal, we have
$F_\H(\cab_{\lambda^*}(\cut(L,L_i)))=x.y_{\lambda^i}$ for some scalar $x\in
R$.  Thus, using the convention given at the beginning of Section \ref{S:MI}
we have $x=\quo{F_\H(\cab_{\lambda^*}(\cut(L,L_i)))}{y_{\lambda^i}}$.
\begin{defi}
  For any $k$-component framed oriented colored link $L$ whose coloring is
  given by a $k$-tuple of $R$-admissible Young diagrams $\lambda^*$ we define
 \begin{align*} H(L,\lambda^*)&=F_\H( \cab_{\lambda^*}(L)), &
 H'(L,\lambda^*,L_i)&=\quo{F_\H\left(
    \cab_{\lambda^*}(\cut(L,L_i))\right)}{y_{\lambda^i}}. 
    \end{align*}
   We call $H$ the colored HOMFLY-PT polynomial and $H'$ the reduced colored
    HOMFLY-PT polynomial.
\end{defi}
Also remark that
$H(L,\lambda^*)=H'(L,\lambda^*,L_i).H(\text{unknot},\lambda^i)$.

\begin{rk}\label{R:framH}
  \begin{enumerate}
  \item \label{RI:framH1} As shown in \cite{AM} the twist acts on
    $y_{\lambda}$ by the scalar
    $$\theta_{\lambda}=a^{|\lambda|^2}v^{-|\lambda|}s^{2n(\lambda)}$$
    where $n(\lambda)=\sum_{(i,j)\in\lambda}j-i=\sum i(\lambda_i-\lambda'_i)$
    (the first sum runs over the coordinates of the cells in the young diagram
    $\lambda$).  So one can re-normalize the colored HOMFLY-PT polynomial to
    get an invariant of colored oriented link (not framed).  For example, if
    all the component of a framed link are colored with the same $\lambda$
    then let $w$
     be the total linking number of $L$ (the algebraic number of crossing of
    $L$) and so, $\theta_{\lambda}^{-w}H(L,\lambda^*)$ is an invariant of the
    underlying oriented link.
    \\
    We can define similarly with the same correction an unframed version of
    $H'$.  Remark that these two unframed invariants are computed with framing
    degree $0$ elements of $\H$ and so they are elements of $\QQ(v,s)$.
  \item In the case of a knot (that is if there is only one component), it
    follows from \cite{M} that $H'(L,\lambda,L)$ is always a Laurent
    polynomial in the variables $a,v,s$ but this is not true for links with
    several components, even for the Hopf link (see \cite{LM}).
  \end{enumerate}
\end{rk}

\begin{prop} \label{P:H2M} For any framed oriented link $L$ with $k$
  components and any $k$-tuple $\lambda^*=(\lambda^1,\ldots,\lambda^k)$ of
  $R$-admissible Young diagrams, one has
\begin{enumerate}
\item $\psi_{m-n}(H(L,\lambda^*))=\quo{\Fmn(L;\qum
    V^{\lambda^*})}{\Id_{\CC[[h]]}}$ \label{PI:H1}
\item $\psi_{m-n}(H'(L,\lambda^*,L_i)) = \quo{\Fmn\left(\cut(L,L_i); \qum
      V^{\lambda^*}\right)}{\Fmn\left(y_{\lambda^i}\right)}$ \label{PI:H2}
\item \label{it:H2M}So, if $\Vschur{\lambda^i}$ is typical then
  $$\psi_{m-n}\left(\frac{H(L,\lambda^*)}{H(\text{unknot},\lambda^i)}\right)
  = \psi_{m-n}\left(H'(L,\lambda^*,L_i)\right) =
  \frac{F'_{m|n}((L;\qum{\Vschur{\lambda^*}}))}
  {\qd(\qum{\Vschur{\lambda^i}})}$$\label{PI:H3}
\end{enumerate}
where $(L;\qum{\Vschur{\lambda^*}})$ denote the link $L\in \T_\Uh$ with its
$i$\textsuperscript{th} component colored by $\qum{\Vschur{\lambda^i}}$, for all
$1\leq i \leq k$.
\end{prop}
\begin{proof}
  The proof of \eqref{PI:H1} follows from the 
  commutative diagram in Proposition \ref{P:H2sl}.  For \eqref{PI:H2} we have
  $F_\H(\cab_{\lambda^*}(\cut(L,L_i)))=x.y_{\lambda^i}$ for some $x \in R$.
  Proposition \ref{P:H2sl} implies that $\wb
  \Fmn(x.y_{\lambda^i})=\psi_{m-n}(x)\wb\Fmn(y_{\lambda^i})$.  Now Proposition
  \ref{P:part2weight} states that $\wb \Fmn(y_{\lambda})$ is a minimal
  projection on an irreducible $\Uh$-module $\qum{\Vschur{\lambda}}\subset\qum
  V^{\otimes r}$.  Thus, the 
  commutative diagram of Proposition \ref{P:H2sl}
  implies \eqref{PI:H2}.  Finally, \eqref{PI:H3} follows from \eqref{PI:H2}
  and the definition of $F'_{m|n}$, since all three quantities in the
  statement correspond to the scalar associated to the (1,1)-tangle whose
  closure is $L$.
\end{proof}
Remark that the first two statements of the proposition are still valid when
$n=0$, i.e. when $\g$ is the Lie algebra $\sll(m)$.
\begin{coro} Let $L$ be an oriented link with $k$ ordered components.  The
  multivariable invariant $M^c_{\sll(m|n)}(L)$ associated with $\sll(m|n)$ and
  $c\in\NN^{m+n-2}$ is determined by the infinite family of colored HOMFLY-PT
  polynomials of $L$.
\end{coro}
\begin{proof}
  First remark that for fixed $c$ there is finitely many $\alpha$ such that
  $\qum V_{\wta^c_\alpha}$ is an atypical module.  So there exists a $
  K_c\in\NN$ such that, for all $ a\in\NN,$ with $a\geq K_c$ we have $\qum
  V_{\wta^c_a}$ is typical.  Now, for $i=1,...,k$, let $a_i\in \NN$ such
  that $a_i\geq K_c$.  Then $M^c_{\sll(m|n)}(L)(q,q^{a_1},\ldots,q^{a_k})$ is
  determined by the linking matrix of $L$ and by $F'_{m|n}(L;(\qum
  V_{\wta^c_{a_1}},\ldots,\qum V_{\wta^c_{a_k}}))$ as all the modules
  $\qum V_{\wta^c_{a_i}}$ are typical.

  Since $M^c_{\sll(m|n)}(L)(q,q_1,\ldots,q_k)$
  (resp. $M^c_0(q,q_1)M^c_{\sll(m|n)}(L)(q,q_1)$ if $k=1$, where
  $M^c_0(q,q_1)$ is defined in the Appendix) is a Laurent polynomial, it is
  determined by the infinite family of polynomials
  $M^c_{\sll(m|n)}(L)(q,q^{a_1},\ldots,q^{a_k})$ for $a_i>K_c$.  Hence to
  prove the corollary it suffices to show that when $a_i\geq K_c$ the
  invariant $F'_{m|n}(L;(\qum V_{\wta^c_{a_1}},\ldots,\qum
  V_{\wta^c_{a_k}}))$ can be computed from the colored HOMFLY-PT polynomial
  of $L$.
  
  To prove this let $a\in\NN,$ with $a\geq K_c$ and choose $(l',l)\in\NN^2$
  such that $m+n+l'-(m+l+\sum_{j=m}^{m+n-2}c_j)=a$, where
  $c=(c_1,...,c_{m+n-2})$.  Let $\mu_1=(\sum_{j=i}^{m-1}c_j)_{i=1\cdots m}$
  and $\mu_2=(\sum_{j=m-1+i}^{m+n-2}c_j)_{i=1\cdots n}$.  We construct the
  partition
  $$\lambda=\left[\left(\underset{m}{\underbrace{n+l',\ldots,n+l'}}\right)
    + \mu_1\right] \cup \underset{l}{\underbrace{n,\ldots,n}} \cup \mu_2'.$$
  where $\mu'_2$ is the conjugate partition of $\mu_2$
  $$\text{i.e.}\qquad\lambda=\pspic{3ex}
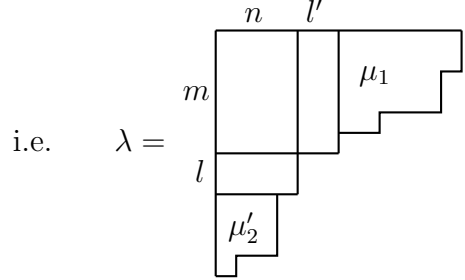
{
    \begin{pspicture}(-1,0)(6,6.5) \psline(0,6)(6,6) \psline(0,6)(0,0)
      \psline(0,3)(3,3) \psline(3,6)(3,3) \psline(0,2)(2,2) \psline(2,6)(2,2)
      \psline(0,0)(0.5,0)(0.5,0.5)(1.5,0.5)(1.5,2)
      \psline(6,6)(6,5)(5.5,5)(5.5,4)(4,4)(4,3.5)(3,3.5) \put(-0.5,2.3){$l$}
      \put(-.8,4.3){$m$} \put(.7,6.2){$n$} \put(2.2,6.2){$l'$}
      \put(0.3,1){$\mu'_2$} \put(3.5,4.8){$\mu_1$}
    \end{pspicture}}$$
  Then one has $\lambda_i-\lambda_{i+1}=c_i$ for $i=1,\dots,m-1$ and
  $\lambda'_i-\lambda'_{i+1}=c_{m-1+i}$ for $i=1,\dots,n-1$.  Furthermore,
  $\lambda_m=n+l'\geq n$ and $\lambda_{m+1}\leq n$.  Thus the module
  $\qum{\Vschur{\lambda}}$ is isomorphic to $\qum V_{\wta^c_a}$.  This with
  \eqref{it:H2M} of Proposition \ref{P:H2M} completes the proof.
\end{proof}

\section{Rank-Level duality in $\H$ and Kashaev's invariant}\label{S:RLK}

Throughout this section, $N$ is an integer greater than $2$ and
$\xi=\e^{\sqrt{-1}\pi/N}$ is the primitive $2N$ root of $1$.  In \cite{ADO} an
invariant of links $L$ whose components are colored by complex numbers is
defined.  This invariant come from the Reshetikhin-Turaev functor of tangles
colored by $U_q(\sll(2))$-finite dimensional modules when $q$ is a root of
unity (see \cite{JM,GPT}).  More precisely, each component of $L$ shall be
colored with the nilpotent representation of $U_q(\sll(2))$ with highest
weight the complex number associated to it.
\begin{defi}
Let $ADO_N$ be the ordered oriented link invariant defined for a
link $L$ with $k$ components and linking matrix $lk_{ij}$ by 
$$
ADO_N(L)(a_1,\ldots,a_k)=q^{-\sum
  lk_{ij}a_i(a_j+2-2N)/2}\Phi^N_L(a_1,\ldots,a_k)
$$
where $q=\e^{\sqrt{-1}\pi/N}$, $a_i\in\CC$ and $\Phi^N_L$ is the framed
ordered link invariant given in \cite{JM}.
\end{defi}
The invariant $ADO_N$ is a slightly modified version of the analogous
invariant constructed in \cite{ADO}.  Contrary to \cite{ADO,MM}, here we
require that the two dimensional representation of $U_q(sl_2)$ has weight $1$.

\begin{theo}[H. and J. Murakami, \cite{MM}, Theorems 4.9 and
  2.1]\label{th:ADO2K}
  For all $ N\geq 2$ and for any link $L$
  $$K_N(L)=J_N(L)_{|q=\xi}=
  ADO_N(L)(\mathsmall{{N-1},\ldots,{N-1}}).$$
\end{theo}
The colored Jones polynomial can be computed as the $\sll(2)$ specialization
of the colored HOMFLY-PT polynomial.  The representation of $\sll(2)$ used by
Murakami is the $N-1$\textsuperscript{th} symmetric power of the standard
representation of $\sll(2)$, i.e. the representration $\Vschur{[N-1]}$ where
$[N-1]$ stands for the partition of $N-1$ consisting in only one part (or the
young diagram with only one row) and $V$ is the standard representation of
dimension $2$. Hence, Theorem~\ref{th:ADO2K} can be restated as follows.
\begin{prop}\label{P:KashH} Let $L$ be an ordered oriented link with $k$
  components.  As before we fix a framed representant of $L$ and we still
  denote it by $L$.  Then
  $$\left.\psi_2\big((\theta_{\scriptscriptstyle[N-1]})^{-w}\,H'(L;
    (\underset{k}{\underbrace{\scriptstyle{[N-1],\ldots,[N-1]}}})
    ;L_i)\big)\right|_{q=\xi} =K_N(L)$$ where $\scriptstyle{[N-1]}$ stands for
  the young diagram corresponding to the one part partition of $N-1$, $w$ is
  the sum of all the linking numbers of $L$ and
  $(\theta_{\scriptscriptstyle[N-1]})^{-w}$ is the framing correction as given
  in Remark \ref{R:framH} \eqref{RI:framH1}, i.e.
  $\theta_{\scriptscriptstyle[N-1]}=a^{(N-1)^2}v^{-(N-1)}s^{(N-1)(N-2)}$.
\end{prop}
\begin{proof}
  We can apply statement \eqref{PI:H2} of Proposition \ref{P:H2M} to
  $\sll(2)$.  As we said, the $N-1$\textsuperscript{th} symmetric power of the
  standard representation of $\sll(2)$ is the $N$-dimensional irreducible
  $\sll(2)$-module whose deformation is used to compute $J_N(L)$.  We just
  need to make the standard correction of $H'$ so that it is a link invariant
  (i.e. framing independent).
\end{proof}

We now introduce a symmetry of $\H$ that leads to the rank-level duality of
modular categories of type $A$.
\begin{defi}
  Let $\Theta: R\go R$ be the ring involution defined by
  $$\Theta(s)=s^{-1} \quad \Theta(a)=-a \quad \Theta(v)=-v$$ 
  We extend $\Theta$ to an involutive functor $\Theta:\T\go\T$ that fixes the
  tangles.
\end{defi}

\begin{prop}\label{P:dual}
  The involution $\Theta$ induces a well defined involutive $\QQ$-linear
  endo-functor of $\H$ such that $\Theta(y_{\lambda})$ is a minimal projector
  of type $\lambda'$ and thus for any framed link with $k$ components,
  $$\Theta(H(L,\lambda^*))=H(L,{\lambda^*}')\quad\text{and}\quad
  \Theta(H'(L,\lambda^*,L_i))=H'(L,{\lambda^*}',L_i)$$ where ${\lambda^*}'$ is
  the $k$-tuple of the conjugate young diagrams of $\lambda^*$.
\end{prop}
\begin{proof}
  First, $\Theta$ is well defined on $\H$ because $\Theta$ fixes the HOMFLY-PT
  skein relations.  In \cite{AM,CB} the symetrizer and antisymetriser of $H_r$
  are defined respectively by:
  $$\,\quad f_r=\frac1{\qnc
    r_s!}s^{-r(r-1)/2}\sum_{\pi\in\Sy(r)}(a/s)^{-l(\pi)}w_\pi$$
  $$\text{and}\quad g_r=\frac1{\qnc
    r_s!}s^{r(r-1)/2}\sum_{\pi\in\Sy(r)}(-as)^{-l(\pi)}w_\pi$$
  where $w_\pi$ are braids of length $l(\pi)$ (the positive permutation
  braids).  Clearly $\Theta(f_r)=g_r$ and $\Theta(g_r)=f_r$.  Now if
  $\lambda=(\lambda_1=l,\lambda_2,\ldots,\lambda_k)$ is $R$-admissible, then
  so is its conjugate partition
  $\lambda'=(\lambda'_1=k,\lambda'_2,\ldots,\lambda'_l)$ (because
  $\qnc{hl(\lambda)}_s=\qnc{hl(\lambda')}_s$).  Furthermore, $c_{\lambda} H_r$
  is the intersection of the two two-side ideals of $H_r$ generated by
  $f_{\lambda_1}\otimes \cdots\otimes f_{\lambda_k}$ and
  $g_{\lambda'_1}\otimes \cdots\otimes g_{\lambda'_l}$.  Hence
  $\Theta(c_{\lambda} H_r)=c_{\lambda'}H_r$ and so $\Theta(y_{\lambda})$ is a
  minimal idempotent of type $\lambda'$.
\end{proof}
\begin{lem}\label{L:Lpsi} For any element $x\in\QQ(s)[v^\pm]$ without a pole
  at $s=\xi$, we have
 \begin{equation}\label{E:Lpsi}
 \left.\psi_\delta\circ\Theta(x)\right|_{q=\xi}
  =\left.\psi_{N-\delta}(x)\right|_{q=\bar\xi}
  \end{equation}
  where $\bar\xi$ is complex conjugation of $\xi$.
  Consequently, if the $k$-tuple of young diagrams $\lambda^*$ satisfy
  $|\lambda^i|<N$ for $i=1\cdots k$, then for a framed link with $k$
  components colored by $\lambda^*$, one has
  $$\left.\psi_\delta((a^{-1}v)^{w}
    H'(L,\lambda^*,L_i))\right|_{q=\xi}= 
  \left.\psi_{N-\delta}((a^{-1}v)^{w}
    H'(L,{\lambda'}^*,L_i))\right|_{q=\bar\xi}$$
  where $w={^t|\lambda^*|}.lk(L).|\lambda^*|=\fdeg(H'(L,\lambda^*,L_i))$
\end{lem}
\begin{proof}
  Equation \eqref{E:Lpsi} follows from a direct calculation.  If all young
  diagrams have size less than $N$, there will be no pole at $s=\xi$ in the
  expression of the associated idempotents $y_{\lambda^i}$.  The factor
  $(a^{-1}v)^{w}$ is invariant by $\Theta$ and $(a^{-1}v)^{w}
  H'(L,\lambda^*,L_i)\in\QQ(s)[v^{\pm1}]$.  Then the last statement of the
  lemma follows from Equation \eqref{E:Lpsi}.
\end{proof}

\begin{lem} \label{L:dKash} Let $V$ be the standard representation of
  $\sll(N-1|1)$  
  and $[1^{,N-1}]=[1,1,\ldots,1]$ be the $N-1$ parts partition of $N-1$ (or
  the young diagram with only one column).  Then
  $\Vschur{[1^{,N-1}]}=\Lambda^{N-1}V=V^0_1$ is the typical representation of
  $\sll(N-1|1)$ with highest weight $(0,\ldots,0,1)$.  Its modified quantum
  dimension $\qd$ has the following properties
  $$\qd_{\sll(N-1|1)}(\qum V^0_1)
  =\prod_{i=1}^{N-1}(q^i-q^{-i})^{-1}= \qn{N-1}!^{-1}$$
  $$\left.\qd_{\sll(N-1|1)}(\qum V^0_1)\right|_{q=\xi}
  =\frac{\e^{-\sqrt{-1}(N-1)\pi/2}}N.$$
\end{lem}
\begin{proof}
  The first statement follows from Proposition \ref{P:part2weight}.  
  The formulas for $\qd$ are taken from \cite{GP2}.
\end{proof}
\begin{theo}\label{th:LG2K} Let $L$ be an oriented link.
  Let $M^0_{\sll(N-1|1)}$ as defined in \cite{GP2} and $\xi=\e^{\sqrt{-1}\pi/N}$.
  Then
  $$K_N(L) = N\e^{\sqrt{-1}(N-1)\pi/2}
  M^0_{\sll(N-1|1)}(L)(\bar\xi,\bar\xi,\ldots,\bar\xi) =
  LG^{N-1|1}(L)(\xi,\bar\xi).$$
\end{theo} 
\begin{proof} Let $L$ be an ordered link.  Fix a framed representant $L_f$ of
  $L$.  We consider the following three colorings of $L_f$.  Let $L'$ be $L_f$
  with all its components colored by the module $\qum V_{\wta^0_1}$.  Let
  $L''$ (resp. $L'''$) be $L_f$ colored by the young diagram corresponding to
  the $N-1$ parts (resp. to the one part) partition of $N-1$.
  Let $w$ is the sum of all the linking numbers of $L_f$.  Let
  $\theta_{[1^{,N-1}]}$ and $\theta_{[N-1]}$ be the values of the twist for
  the colored HOMFLY-PT given in Remark \ref{R:framH} \eqref{RI:framH1} and
  let $\theta^0_1$ be the value of the twist on $\qum V_{\wta^0_1}$.
  
  With this notation we have 
  $$\begin{array}{rl}
    N\e^{\sqrt{-1}(N-1)\pi/2}
    M^0_{\sll(N-1|1)}(L)(\bar\xi,\bar\xi,\ldots,\bar\xi)=& 
    \left.\qd_{\sll(N-1|1)}(\qum V_{\wta^0_1})^{-1}
      (\theta^1_0)^{-w}F'(L')\right|_{q=\bar\xi}
    \\=&\left.\psi_{N-2}\left(\theta_{[1^{,N-1}]})^{-w}
        H'(L'')\right)\right|_{q=\bar\xi} 
    \\=&\left.\psi_2\left(\Theta((\theta_{[1^{,N-1}]})^{-w}
        H'(L'')\right)\right|_{q=\xi} 
    \\=&\left.\psi_2\left((\theta_{[N-1]})^{-w}H'(L'''\right))\right|_{q=\xi}
    \\=&K_N(L)
  \end{array}$$
  where the first equality follows from the definition of $M^0_{\sll(N-1|1)}$
  and Lemma \ref{L:dKash}, the second from Proposition \ref{P:H2M}
  (\ref{it:H2M}), the third from Lemma \ref{L:Lpsi}, the fourth from
  Proposition \ref{P:dual} and the last equality follows from Proposition
  \ref{P:KashH}. 
\end{proof}
We are lead to the following conjecture.
\begin{conj}\label{cj:M2ADO}
  Let $\xi=\e^{\sqrt{-1}\pi/N}$ then
  $M^0_{\sll(N-1|1)}(L)(\xi,\xi^{a_1},\ldots,\xi^{a_c})= ADO_N(a_1,\ldots,a_c).$
\end{conj}
\section*{Appendix}
In this appendix we set $\g=\sll(m,1)$.  In Lemma 3.2 of \cite{GP2} we define
two Laurent polynomials $M_0^c\in \ZZ[q^{\pm1}]$ and $M_1^c\in \ZZ[q^{\pm1},
q_1^{\pm1}]$ and show that
$$\qd(V_{\wta^c_{\alpha}})=\frac{M_0^c(e^{h/2})}{M_1^c(e^{h/2}, e^{\alpha h/2})}$$ 
where $c\in \NN^{m-1}, \alpha \in \CC$.  When $c=0$ we have
$M_1^0(q,q_1)=\prod_{i=0}^{m-1}(q_1q^i-(q_1q^i)^{-1})$ and $M_0^0(q)=1$.
Therefore, we have the following lemma.
\begin{lem}\label{L:d(V_a)1} 
  $\qd(\qum V_{\wta^0_{\alpha}})=\frac{1}{\prod_{i=0}^{m-1}
    (q^{\alpha+i}-q^{-(\alpha+i)})}.$
\end{lem}

Now a direct calculation gives the following lemma.
\begin{lem}\label{L:d(V_a)}
  $M_1^0(\e^{\sqrt{-1}\pi/m},q_1)=\e^{\sqrt{-1}(m-1)\pi/2}(q_1^m-q_1^{-m}).$
\end{lem}

Let $V_l$ be the $\g$-modules determined by Equation \ref{Eq:split}.  For
$0<l<m$ we have $V_l=V_{\wta^{c_l}_{2\alpha+m-i-1}}$ where $c_l$ is the
$l$-tuple with a $1$ in the $l$ slot and $0$ everywhere else.
\begin{lem}\label{L:d(V_i)} When $q=e^{\sqrt{-1}\pi/m}$ we have $\qd(V_l)=0$
  for $0<l<m$.
\end{lem}
\begin{proof}
Using the notation and results in the Appendix of \cite{GP2} we have 
$$\rho=\sum (m-i)\epsilon_i,\;\; <\rho,\epsilon_i-\epsilon_j>=j-i, \;\;
<w_k,\epsilon_i-\epsilon_j>=\left\{\begin{array}{l} 1 \;\; \text{ if } i\leq
    k<j,\\ 0 \;\; \text{ else.}\end{array}\right.$$ These equalities and the
fact that the Laurent polynomials $M_0^c$ and $M_1^c$ are equal to particular
specializations of elements in the group ring of the weight lattice (see Lemma
3.2 of \cite{GP2}) imply
$$M_0^{c_l}(q)=\dfrac{\qn{m}!}{\qn{l}!\qn{m-l}!}= 
\dfrac{\prod_{i=1}^l\prod_{j=l+1}^m\qn{j-i+1}}
{\prod_{i=1}^l\prod_{j=l+1}^m\qn{j-i}}.$$ Therefore, when
$q=e^{\sqrt{-1}\pi/m}$ we have
$\qd(V_l)=M_0^{c_l}(q)/M_1^{c_l}(q,q^{2\alpha+m-l-1})=0$.
\end{proof}

\linespread{1}

\vfill


\begin{thebibliography}{99}
\setcounter{bibcount}{0}

\bibitem{AM} A. K. Aiston, H. R. Morton - {\em Idempotents of Hecke algebras
    of type $A$}, J. Knot Theory Ramifications {\bf 7} (1998), no.~4,
  463--487.

\bibitem{ADO} Y. Akutsu, T. Deguchi, and T. Ohtsuki - {\em Invariants of
    colored links.} J. Knot Theory Ramifications \textbf{1} (1992), no. 2,
  161-184.

\bibitem{CB} C. Blanchet - {\em Hecke algebras, modular categories and
    $3$-manifolds quantum invariants.} Topology \textbf{39} (2000), no. 1,
  193--223.

\bibitem{WIL} D. De Wit, A. Ishii, J.R. Links - {\em Infinitely many
    two-variable generalisations of the Alexander-Conway polynomial.}  Algebr.
  Geom. Topol.  5 (2005), 405--418.

\bibitem{D} T. Deguchi - {\em Multivariable invariants of colored links
    generalizing the Alexander polynomial.} Proceedings of the Conference on
  Quantum Topology (Manhattan, KS, 1993) (River Edge, NJ), World Sci.
  Publishing, 1994, 67--85.

\bibitem{G04A} N. Geer - {\em Etingof-Kazhdan quantization of Lie
    superbialgebras.}  { Advances Math.}  
  \textbf{207} (2006), no. 1, 1--38.

  
\bibitem{GP1} N. Geer, B.Patureau-Mirand - {\em Multivariable link invariants
    arising from $\sll(2|1)$ and the Alexander polynomial.}  J. Pure
  Appl. Algebra 210 (2007), no. 1, 283--298.  (math.GT/0601291).

\bibitem{GP2} N. Geer, B.Patureau-Mirand - {\em Multivariable link invariants
    arising from Lie superalgebras of type I.}  Preprint (math.GT/0609034).

\bibitem{GPT} N. Geer, B.Patureau-Mirand, V. Turaev - {\em Modified quantum
    dimensions and re-normalized link invariants.}  Preprint.

\bibitem{K} V.G. Kac - {\em Lie superalgebras.} Advances Math. 26 (1977),
  8--96.
  
\bibitem{Kv} R. M. Kashaev - {\em A link invariant from quantum dilogarithm.}
  Modern Phys. Lett. A \textbf{10} (1995), no. 19, 1409--1418.
  
\bibitem{LG} J.R. Links, M. Gould - {\em Two variable link polynomials from
    quantum supergroups.} Lett. Math. Phys. \textbf{26} (1992), no. 3,
  187--198.

\bibitem{LM} S.G. Lukac, H.R. Morton - {\em The HOMFLY polynomial of the
    Decorated Hopf Link} J. Knot Theory Ramif. 12 (2003) 395-416.

\bibitem{M} H.R. Morton - {\em Integrality of Homfly (1,1)-tangle invariants}
  Algebraic and Geometric Topology, 7 (2007), 227-238.

\bibitem{MM} H. Murakami, J. Murakami - {\em The colored Jones polynomials and
    the simplicial volume of a knot.}  Acta Math. 186 (2001), no. 1, 85--104.

\bibitem{JM} J. Murakami - {\em Colored Alexander Invariants and
    Cone-Manifolds.}  Acta Math. 186 (2001), no. 1, 85--104. 

\bibitem{Se} A.N. Sergeev - {\em Tensor algebra of the identity representation
    as a module over the Lie superalgebras ${\rm Gl}(n,\,m)$ and $Q(n)$.} Mat.
  Sb. (N.S.) 123(165) (1984), no. 3, 422--430.

\bibitem{RT}N.Y. Reshetikhin, V.G. Turaev - {\em Ribbon graphs and their
    invariants derived from quantum groups}, Comm. Math. Phys. {\bf 127}
  (1990), no.~1, 1--26.

\bibitem{Tu86} V.G. Turaev - {\em Reidemeister torsion in knot theory.}  [J]
  Russ. Math. Surv. 41, No.1, 119-182 (1986). (English translation of Uspekhi
  Mat. Nauk 41 (1986), no. 1(247), 97--147, 240.)
  
\bibitem{Tu94} V.G. Turaev - {\em Quantum invariants of knots and
    3-manifolds.} de Gruyter Studies in Mathematics, 18. Walter de Gruyter \&
  Co., Berlin, (1994) 588 pp.
  
\end{thebibliography}
\end{document}